
\documentstyle{amsppt}


\nopagenumbers

\input amstex
\magnification=\magstep1

\define\sbh{\subheading}

\define\fm{{\Im}_M}

\redefine\LL{{\Cal L}}

\define\qc{\text{quasi-coherent}}
\define\Aq{A_{qc}-mod}
\define\cxa{C_X(A)}
\define\cxan{C_X^{new}(A)}
\define\cxao{C_X^{old}(A)}\define\shx{Sh_X(A)}
\define\shxn{Sh_X^{new}(A)}
\define\shxo{Sh_X^{old}(A)}
\define\qcsh{Sh_X^{qc}(A)}
\define\fs{f_{*}}
\define\fd{f^{\bullet}}
\redefine\no{new\to old}
\redefine\on{old\to new}
\define\Dq{D_{qc}(A)}
\define\Dmq{D^-_{qc}(A)}
\define\Ds{D(Sh_X^{qc}(A))}
\define\Dms{D^-(Sh_X^{qc}(A))}
\define\Dsq{D_{qc}Sh_X(A)}
\define\Dbs{D^b(Sh_X^{qc}(A))}

\define\Dbsq{D^b_{qc}Sh_X(A)}

\topmatter
\title Grothendieck topologies \\
and deformation theory II
\endtitle
\author  Dennis Gaitsgory \endauthor
\address{School of Mathematical Sciences, Tel-Aviv University,
Ramat-Aviv, Israel} \endaddress
\email
gaitsgde\@math.tau.ac.il
\endemail
\endtopmatter
\heading 0. Introduction \endheading
\sbh{0.1}
In the present paper we continue the study of deformation theory
of algebras using the approach of [Ga]. We will extend the main results
of [Ga] to the global case. Namely, we pose and solve the following
problem: what cohomological machinery controls deformations of a sheaf of
algebras over a scheme? This question has already been studied by many
authors [Ill],[Ge],[H-Sch],[Schl].
\sbh{0.2}
Let first $A$ be an associative algebra over a ring. Consider the category
of all algebras over $A$, let us call it $C(A)$. One can observe that every
question concerning the deformation theory of $A$ can be formulated in terms
of this category.

Our first step will be to apply a linearization
procedure to $C(A)$, in other words we will endow it with a Grothendieck
topology and then we will consider sheaves of abelian groups on it. It
will turn out that deformations of $A$ are controlled by cohomologies
of certain sheaves on this site. Cohomologies arise naturally as classes
attached to torsors and gerbes. All this was done in [Ga].

When $A$ is
no longer an algebra over a ring but rather a $\qc$ sheaf of algebras over
a scheme $X$, the definition of $C(A)$ must be modified in order to take
into account possible localization with respect to $X$, since the
appropriate cohomology theory would incorporate algebra cohomology
of $A$ and scheme cohomology of $X$. In this case
instead of working with the whole category of sheaves on our site, we
single out a subcategory which we call the category of $\qc$ sheaves.
This category will have properties similar to those of the category
of $\qc$ sheaves of $A-$bimodules among all sheaves of $A-$bimodules and
it will be more manageable.

The second step will be to find a connection between the category
of sheaves on $C(A)$ and the category of $\qc$ sheaves of $A-$bimodules on
our scheme $X$. This connection will be described by two mutually
adjoint functors, which would enable us to rewrite the cohomology
groups that control deformations of $A$ in terms of cohomologies of
some canonical object $T^{\bullet}(A)$ of the derived category of
$\qc$ sheaves of $A-$bimodules. The object $T^{\bullet}(A)$ will be called
the cotangent complex of $A$. Another approach to the construction of
the cotangent complex in a slightly different situation was used by Illusie
[Ill].

\sbh{0.3} Let us now describe the contents of the paper.

In Section 1 we present a brief exposition of some well known
facts and results from the theory of sites. For a more detailed
discussion the reader is referred to [Ar,Gr]. In the remaining sections
we will freely operate with the machinery of sheaves, cohomologies,
direct and inverse images; therefore the reader is advised to look
through this section in order to become familiar with the notation.

In Section 2 we define the site $\cxa$ along with its variants for
affine schemes. We introduce also the appropriate categories of sheaves
and functors between them. The central results are
\roster
\item
Theorem 2.3.3  with its corollaries, that insure that the category
$Sh^{qc}(A)$ is well defined
\item
Theorem 2.5 that says that cohomologies of $\qc$ sheaves computed inside
the $\qc$ category and inside the category of all sheaves give the
same answers.
\endroster

In Section 3 we introduce functors $\Im$ and $\LL$ that establish
connection between the category $Sh^{qc}(X)$ and the category $\Aq$.
Let us remark that it would be possible to work with the category of
all sheaves on $\cxa$ without introducing $\qc$ sheaves explicitly.
We, nevertheless, decided to that, since to our mind, introducing
this category and basic functors that are connected to it reflects
the nature of the things and clarifies the exposition.

Finally, Section 4 is devoted to deformation theory. Theorem 4.1.2
describes how to pass from deformations to cohomology of sheaves
on $\cxa$ via torsors and gerbes, and in 4.2 we translate the
assertions of this theorem to the language of cohomology of $\qc$
sheaves of $A$-bimodules.

\sbh{0.4}
The results of the present paper can be easily generalized to the
case of algebras over an arbitrary operad (cf. [Ga,H-Sch]). We opted for
treating the case of associative algebras only in order to simplify the
exposition.
One can also develop a similar theory for operad co-algebras.

\sbh{0.5}
In recent years there have been a lot of interest in deformation
theory. We have to mention the works [H-Sch,Ge-Sch,Ma,St-Schl,Fo]. Our
approaches connected to that of [Ill].
Let us also point out that one of the
central ideas of the present paper: to resolve an algebra $A$ by
free algebras (at least locally) goes back probably to Quillen
and to Grothendieck [Qu].
\sbh{Acknowledgments}
The author is deeply indebted to J.Bernstein under whose supervision
the work was written.
Also I would like to thank L.Breen, A.Joseph, and S.Shnider for interesting
and stimulating discussions.

\heading {1. preliminaries on Grothendieck topologies} \endheading
\sbh{1.0}
In this section we will review certain notions from the theory
of sites. Proofs will be given mostly in cases when our exposition
differs from the standard one.
\sbh{1.1}
Let $C$ be a category possessing fiber products.
A Grothendieck topology (cf. [Gr]) on it (or a structure of a site)
is a collection of morphisms that are called covering maps
if it satisfies the following three conditions:
\roster
\item
Any isomorphism is a covering.
\item
If $\phi:U\to V$ and $\psi:V\to W$ are coverings, then their composition
$\psi\circ\phi: U\to W$ is a covering too.
\item
If $\phi:U\to V$ is a covering and if $\alpha:V_1\to V$ is an arbitrary
morphism, then the base change map $\phi_1:U\underset V\to\times V_1\to V_1$
is a covering.
\endroster
\sbh{1.1.1 Examples}
\sbh{1} For any category $C$ there exists the minimal Grothendieck topology:
the only coverings are isomorphisms. This site will be denoted by $(C,min)$.

\sbh{2} Let $Set$ be the category of sets. We introduce the structure of
a site on it by declaring surjections to be the covering maps.

\sbh{3} Let $Set^o$ be the category opposite to $Set$. We introduce a
Grothendieck
topology by declaring $\phi:X\to Y$ to be a covering
if the corresponding map of sets $Y\to X$ is an injection.

\sbh{4} Constructions similar to the above
ones can be carried out when the category $Set$
is replaced by an abelian category, in particular, by the category
$Ab$, the category of abelian groups.

\sbh{5} Let $X$ be a topological space. Let $C(X)$ be the category whose
objects are
finite disjoint unions of open subspaces of $X$.
$$Hom(U,V)\overset\text{def}\to=
\text{ maps from $U$ to $V$ compatible with an embedding to $X$}.$$
A map $\phi\in Hom(U,V)$ is a covering if it is surjective.

\sbh{6}
Let a cite $C$ have a finite object $X_0$ and let $X$ be any other object
of $C$. We can define a new cite $C_X$ whose underlying category
is the category of "objects of $C$ over $S$", with morphisms being compatible
to the projection to $X$. A morphism $\phi$ in $C_X$ is declared to be
a covering if it is a covering in $C$.

\sbh{1.1.2}
Let $C_1$ and $C_2$ be two sites. A functor $F$ between the underlying
categories is said to be a functor between sites if the following holds:
\roster
\item
$F$ maps coverings to coverings.
\item
If $A,B,D$ are three objects in $C_1$ with $A,B$ mapping to $D$, then
the canonical map $F(A\underset D\to\times B)\to
F(A)\underset F(D)\to\times F(B)$ is a covering in $C_2$.
\endroster

We say that a functor $F$ between two sites is strict if it preserves
fiber products, i.e. if the map in (2) is an isomorphism.

\sbh{1.2}
\sbh{Definition}
A sheaf of sets (resp. of abelian groups) on a site $C$
is a functor $S$ between the sites $C$ and $Set^o$ (resp. $Ab^o$), the latter
considered with the topology specified in the Example 3 above.

Morphisms between sheaves are by definition natural transformations
between such functors.

\sbh{Definition}
A presheaf of sets on $C$ (resp. of abelian groups) is sheaf on $C$
when the latter is considered with the minimal topology.

It is an easy exercise to verify that the above definition of a
sheaf coincides with the traditional one. From now on
by a sheaf we will mean a sheaf of abelian groups. It will be left
to the reader to make appropriate modifications for sheaves of sets.

The category of sheaves will be denoted
by $Sh(C)$. This category possesses a natural additive structure and
is in fact an abelian category. If $S$ is
a sheaf, and if $X\in C$, $S(X)$ will be denoted by $\Gamma(X,S)$
and will be called the set of sections of $S$ over $X$. The map
$\Gamma(X,S)\to\Gamma(Y,S)$ for a map $Y\to X$ will be called
the restriction map.

\sbh{1.3}
Let $F:C_1\to C_2$ be a functor between sites.
We have then the natural functor (called direct image)
$F_{\bullet}:Sh(C_2)\to Sh(C_1)$. This functor is always left exact. It is also
right exact if the following condition is satisfied:

For any covering $Z\to F(X)$ there exists a covering $\phi:Y\to X$,
endowed with a map $\alpha:F(Y)\to Z$ such that the composition
$F(Y)\to Z\to F(X)$ coincides with $F(\phi)$.

The functor $F_{\bullet}$ has a left adjoint (called the inverse
image): $F^{\bullet}:Sh(C_1)\to Sh(C_2)$. The functor $F^{\bullet}$ is
always right exact and it is also left exact if the functor $F$ is strict
in the sense of 1.1.1.

\sbh{1.3.1 Examples}
\sbh{1} Let $Forget:(C,min)\to C$ be the canonical functor of sites.
The above constructions
yield the embedding functor from sheaves to presheaves and
its left adjoint, which is called the functor of associating a
sheaf to a presheaf. It is a good exercise to describe the associated
sheaf explicitly.

\sbh{2} Let $pt$ be the category of one object and one morphism. If $C$
is a site, for any $X\in C$ we have a functor $pt_X:pt\to C$, that sends
the unique object of $pt$ to $X$. We have the canonical constant sheaf
$Const$ on $pt$. Let by definition $Const_X=pt_X^{\bullet}(Const)$. This sheaf
will
be called the constant sheaf corresponding to $X$. By definition
we have: $Hom(Const_X,S)=\Gamma(X,S)$ functorially with respect to
$S\in Sh(C)$.

\sbh{3} Let $F:C_1\to C_2$ be a functor between sites and let $X\in C_1$.
Then $$F^{\bullet}(Const_X)\simeq Const_{F(X)}.$$

\sbh{4} Recall the situation of 1.1.1 Example 6. We have the natural embedding
functor $i: C_X\to C$ and its right adjoint $Cart:Y\to Y\underset{X_0}\to\times
X$.
We claim then, that the functors $i_{\bullet}$ and $Cart^{\bullet}$ are
canonically isomorphic. We denote this functor by $S\to S|C_X$ and call it
the functor of restriction of a sheaf to $C_X$. By definition, for $Y\in C_X$
we have $\Gamma(Y,S|C_X)\simeq\Gamma(Y,S)$.

If now $X\to X_0$ is a covering, the functor $S\to S|C_X$ is exact and
faithful.

\sbh{1.4 Cohomology of sheaves}
Along with the abelian category $Sh(C)$ one considers also the
corresponding derived categories $D(Sh(C))$, $D^+(Sh(C))$,
$D^-(Sh(C))$ and $D^b(Sh(C))$.
It can be shown [Ar,Gr] that the category $Sh(C)$ has enough injective objects.
In particular, any left exact functor admits a right derived functor.
If $X\in C$, $R^i\Gamma(X,S)$ will be denoted by $H^i(X,S)$.

\sbh{1.4.1 $\check{C}$ech complexes}
Let $C$ and $C_X$ be as in 1.3.1 Example 4 above. Put $U_0=X$ and
let $U_i$ denote the $i+1$-fold fiber product of $X$ with itself
over $X_0$. Let also $j_i$ denote the canonical map $j_i:U_i\to X_0$.

For any sheaf $S\in Sh(C)$ we can form a canonical complex:
$$0\rightarrow S\rightarrow {j_0}_{\bullet} {j_0}^{\bullet}(S)\rightarrow
{j_1}_{\bullet}j_1^{\bullet}(S)\rightarrow {j_2}_{\bullet}j_2^{\bullet}(S)
\rightarrow ...$$

\proclaim{Claim}
This complex is exact.
\endproclaim

To prove this statement, we restrict this complex to $X=U_0$
and this enables us to write down an explicit homotopy operator.

The complex
$$0\rightarrow {j_0}_{\bullet} {j_0}^{\bullet}(S)\rightarrow
{j_1}_{\bullet}j_1^{\bullet}(S)\rightarrow {j_2}_{\bullet}j_2^{\bullet}(S)
\rightarrow ...$$
will be called the $\check{C}$ech complex of $S$.

\sbh{1.5 Torsors and Gerbes}
\sbh{1.5.0}
Let now our category possess a final object $X_0$ and let $S$ be
a sheaf of abelian groups. $H^i(S)$ will denote $R^i\Gamma(X_0,S)$.
\sbh{1.5.1}
Before defining torsors and gerbes in the sheaf-theoretic context
we need to recall several definitions.

Let $\Gamma$ be an abelian group and let $\Gamma$ act on a set $\tau$.
We say that $\tau$ is a torsor over $\Gamma$ if this action is simply
transitive. Torsors over a given group form a rigid monoidal category
(cf. [DM])
under $\tau_1\otimes\tau_2\to\tau_1\times\tau_2/\Gamma$ with the
anti-diagonal action of $\Gamma$.

Let now $O$ be a monoidal category and let $M$ be an arbitrary category.
We say that $O$ acts on $M$ if we are given
\roster
\item
A functor $Action:O\times M\to M$.
\item
A natural transformation between the two functors $O\times O\times M\to M$:
$$
\CD
       O\times O\times M    @>{Action}>>   O\times M        \\
          @VVV                             @V{Action}VV     \\
         O\times M           @>{Action}>>    M
\endCD
$$
\endroster
such that the obvious "pentagon" identity is satisfied.

We say that $M$ is a gerbe bound by $O$, if for any $X\in M$
the functor $O\to M$ given by $A\to Action(A\times X)$ is an
equivalence of categories.

If $O$ is a groupoid and if $M$ is a gerbe bound by $O$, then
$M$ is also a groupoid and $\pi_0(M)$ is a torsor over $\pi_0(O)$.
\sbh{1.5.2}
A sheaf of sets $\Upsilon$ is called a torsor over $S$ if
\roster
\item
$S$ viewed as a group-like object in the category of sheaves of
sets acts on the object $\Upsilon$, i.e if for every $X\in C$,
$\Gamma(X,S)$ acts on $\Gamma(X,\Upsilon)$ in a way compatible with
restrictions.
\item
For every $X\in C$, $\Gamma(X,\Upsilon)$ is a torsor over $\Gamma(X,S)$,
whenever the former is nonempty.
\item
For some covering $X$ of $X_0$, the set $\Gamma(X,\Upsilon)$ is nonempty.
\endroster

Let $T_S$ denote the category of torsors over $S$. From 1.5.1
it is easy to deduce that the groupoid $T_S$ possesses a structure
of a rigid monoidal category.
\proclaim{Claim}
The group $\pi_0(T_S)$ is canonically isomorphic to $H^1(S)$.
\endproclaim
\demo{Proof}
In fact, we claim more:
Consider the category $Ext(Const_{X_0},S)$, whose objects are short exact
sequences $0\to S\to E\to Const_{X_0}\to 0$. Then this category is canonically
equivalent to $T_S$.

Indeed, for any such extension $0\to S\to E\to Const_{X_0}\to 0$ we associate
a torsor $\Upsilon$ by setting for every $X$ over $X_0$
$$\Gamma(X,\Upsilon)=\text{ splittings: }\Gamma(X,Const_{X_0})\to\Gamma(X,E)$$

This functor is easily seen to be an equivalence of (monoidal) categories
and \newline ${\pi}_0(Ext(Const_{X_0},S))\simeq H^1(S)$.

\enddemo{QED}
\sbh{1.5.3}
We are heading towards the definition of gerbes, but first
we need to recollect the notion of a stack.

Let $C$ be a site with a final object $X_0$. Suppose that for each $X\in C$ we
are given a category $G(X)$, for each map $\alpha:Y\to X$ we are
given a functor $G_{\alpha}:G(X)\to G(Y)$ and for each composition
of maps $\alpha:Y\to X$ and $\beta:Z\to Y$, we are given a natural
transformation $F_{\alpha}\circ F_{\beta}\to F_{\alpha\circ\beta}$,
such that all the data are compatible with respect to two-fold compositions.
This collection is called a stack if the following two axioms are
satisfied:
\roster
\item
Let $X\in C$, and let us consider the category $C_X$ as in 1.1.1 Example 6.
Let also $s_1,s_2$ be two objects of $G(X)$.
We can consider the presheaf of sets on $C_X$: $Y\in C_X\to Hom(s_1|Y,s_2|y)$.
We require that this presheaf is a sheaf for each $X\in C$.
\item
Let $\phi:Y\to X$ be a covering. Consider the category of descent
data on $Y$ with respect to $X$, whose objects are pairs
$s\in G(Y)$ and an isomorphism $p_1^{*}(s)\to p_2^{*}(s)$, where
$p_1,p_2$ are the two projections from $Y\underset X\to\times Y$
to $Y$, such that the above isomorphism satisfies
the obvious cocycle condition on the three-fold fiber product of
$Y$ with itself over $X$. Morphisms in this category are defined
to be maps $s_1\to s_2$ compatible with isomorphisms between their
pull-backs on $Y\underset X\to\times Y$.
We have the obvious functor from $G(X)$ to this category of descent data.
We require that this functor is an equivalence of categories.
\endroster

\sbh{Examples}
\sbh{1} All sites that we will be working
with in this paper will have the following property:
$G(X):=C_X$ is a stack.

\sbh{2} If $S$ is a sheaf of groups, we can define $G(X)=T_{S|X}$ (torsors
over $S$ restricted to the category $C_X$).
This is a stack too.
\sbh{1.5.4 Gerbes}
Let once again $C$ be a site with a final object $X_0$ and let
$S$ be a sheaf of abelian groups.
Let $G$ be a stack on $C$ endowed with the following additional
structure:
\roster
\item
Each $G(X)$ is acted on by the monoidal category $T_{S|X}$.
\item
For each $\alpha:Y\to X$ we are given a natural transformation between
two functors $T_{S|X}\times G(X)\to G(Y)$:
$$G_{\alpha}\circ Action_X\to Action_Y\circ
({T_{S|X}}_{\alpha}\times G_{\alpha}),$$
which is compatible with the natural transformations of 1.5.1(2) and
with composition of restrictions.
\endroster
Suppose that for each $X\in C$,
$G(X)$ is a gerbe over $T_{S|X}$ and that there exists a covering $X$
of $X_0$ such that $G(X)$ is nonempty. We say then that $G$
is a gerbe bound by $S$.

Functors between gerbes bound by a sheaf of abelian groups $S$
and natural transformations between such functors are defined in
a natural fashion.

\sbh{Remark}
Let $S$ be a sheaf of abelian groups and let $G$ be a stack such
that if $s_1,s_2\in G(X)$, there exists a covering $Y$ of $X$
such that the pull-backs of $s_1$ and $s_2$ on $G(Y)$ become isomorphic.
Then $G$ is a gerbe bound by $S$ if for every $X,s\in G(X)$, $Aut(s)$
is isomorphic to $\Gamma(X,S)$ functorially in $X$ and in $s$.

\sbh{Examples}

\sbh{1}
Let $S$ be in $Sh(C)$. A basic example of a gerbe bound by $S$ is provided
by setting $G(X)=T_{S|X}$. This gerbe of $S$-torsors we denote by abuse of
notation by $T_S$.

It is an easy observation that a gerbe $G$ is equivalent to $T_S$ if
and only if $G(X_0)$ is nonempty.
\sbh{2}
Let $S_1\to S_2$ be a map of sheaves of abelian groups. If $G$ is a gerbe
bound by $S_1$, we can construct an induced gerbe $G'$ bound by $S_2$.
\sbh{3}(cf. [D-Ill,BB])
Let $0\to S\to K_1\to K_2\to Const_{X_0}\to 0$ be an exact sequence of sheaves
on $C$. Let $K^{\bullet}$ denote the 2-complex $K_1\to K_2$. To this
2-complex we can associate a gerbe $G(K_1\to K_2)$ bound by $S$ in a
canonical way by setting: $G(K_1\to K_2)(X)=\text{the category of extensions}$
$0\to K_1|C_X\to E\to Const_X\to 0$ of sheaves over $C_X$ endowed with a map
$E\to K_2|C_X$ commuting with the projection to $Const_{X_0}|C_X\simeq
Const_X$.
It is easy to verify that $G(K_1\to K_2)$ defined in this way is indeed a
gerbe.

If now $\alpha:K^{\bullet}\to {K'}^{\bullet}$ is a quasi-isomorphism of
2-complexes,
we get a canonical functor between the corresponding gerbes $G(K^{\bullet})$
and $G({K'}^{\bullet})$. This means that the operation of assigning a gerbe to
a
2-complex is well defined on the derived category $D(Sh(C))$.

\sbh{1.5.5}
The following proposition is not difficult to prove:
\proclaim{Proposition}
The assignment $K^{\bullet}\to G(K^{\bullet})$ establishes a one-to-one
correspondence between the set isomorphism classes of objects
$K^{\bullet}$ in $D(Sh(C))$ with nontrivial cohomologies only in degrees
$0$ and $1$, such that $H^0(K^{\bullet})\simeq S$,
$H^1(K^{\bullet})\simeq Const_{X_0}$ and the set of equivalence classes of
gerbes $G$ bound by $S$.
\endproclaim

In particular, since the set of isomorphism classes of 2-complexes of the
above type in the derived category is $Ext^2(Const_{X_0},S)\simeq H^2(S)$,
to any gerbe $G$ bound by $S$ we can associate a well defined class in $H^2(S)$
that vanishes if and only if $G(X_0)$ of this gerbe is nonempty.

\heading 2. the site $\cxa$ \endheading
\sbh{2.0}
As it has been explained in the introduction, our bridge between
deformations and cohomology
is based on considering sheaves on the site $\cxa$ which we are
about to define. Throughout this paper, by a scheme we will mean
a separated scheme. It is not difficult, however, to generalize
all our results to the case of an arbitrary scheme.
\sbh{2.1}
Let $X$ be a scheme and let $Zar_X$ denote the Zariski site of $X$,
whose objects are disjoint finite unions of open subsets of $X$ and whose
morphisms
are maps of schemes over $X$. A morphism in $Zar_X$ is a covering
if it is surjective.
Let $A$ be a $\qc$ sheaf of associative
algebras on $X$.
\sbh{2.1.1 Definition of $\cxa$}

Objects: triples $(U,B_U,\phi)$, with $U\in Zar_X$, $B_U$ is a $\qc$
sheaf of associative algebras on $U$ and $\phi:B_U\to A|U$ is a map
of sheaves of associative algebras. Here $A|U$ is the restriction of $A$ on
$U$.
When no confusion can be made, we will omit $\phi$.

Morphisms: $Hom((V,C_V),(U,B_U))$ is a set of pairs $(j,\alpha)$, \newline
where $j\in Hom_{Zar_X}(V,U)$ and $\alpha:C_V\to B_U|V$
(restriction by means of $j$).

The category $\cxa$ is easily seen to have fiber products.

Topology: $(j,\alpha):(V,C_V)\to (B,B_U)$ is said to be a covering map
if $j$ is a covering in $Zar_X$, and if $\alpha$
is an epimorphism.

\smallskip
Sometimes when no confusion can be made we will write $A$ instead of
$\Gamma(X,A)$ for $X$ being an affine scheme.

\sbh{2.1.2 Variant}
When $X$ is an affine scheme $X=Spec(R)$, the site $\cxa$ will often
be denoted by $\cxan$ to emphasize the difference between
$\cxan$ and $\cxao$:
\sbh{ Definition of $\cxao$}

Objects: $R$-algebras $B$ with a map to $\Gamma(X,A)$.

Morphisms: Algebra homomorphisms commuting with \newline
structure maps to $\Gamma(X,A)$.

Topology: Covering maps are defined to be just epimorphisms of algebras.

\sbh{2.1.3}
$\shx$ will denote the category of sheaves of abelian groups over $\cxa$.
For $X$ affine, $X=Spec(R)$, this category will also be denoted by $\shxn$,
whereas $\shxo$ will denote the category of sheaves of abelian groups
over $\cxao$.

\sbh{Example 1}
Let $(U,B_U)\in\cxan$. According to 1.3.1(2) we can consider
the sheaf $Const_{(U,B_U)}\in\shxn$. It follows from the Example 1.5.3(1)
that $Const_{(U,B_U)}$is given by  $(V,C_V)\to\underset{Hom((V,C_V),(U,B_U))}
\to\oplus({\Bbb Z})$.
\sbh{Example 2}
A similar construction can be carried out in the case of
an affine scheme $X=Spec(R)$ for $\cxao$. For a projective $R$-module $V$
let $FreeR(V)$ denote the free associative algebra built on $V$.
If now $Free(V)\in\cxao$, the sheaf $Const_{Free(V)}$ is a projective object
of $Sh^{old}_A(X)$. This is because every covering of $Const_{Free(V)}$
admits a section.

\sbh{2.2}
Let $f:Y\to X$ be a morphism of schemes. Let us be given $\qc$
sheaves of algebras $A$ on $X$ and $A'$ on $Y$. Assume also be given
a map of sheaves of algebras $\phi:A'\to f^{*}(A)$. (Here $f^{*}$ denotes the
pullback.)
We say then that $(f,\phi)$ is a map from the pair $(Y,A')$ to the pair
$(X,A)$.

\sbh{2.2.1} We have a functor denoted $(f,\phi)$ or just $f$:
$$\cxa\to C_Y(A'):(U,B_U)\text{ goes to } (U\underset X\to\times Y,f^{*}B
\underset f^{*}(A|U)\to\times (A'|U\underset X\to\times Y))$$
This functor is strict if $f$ is flat.

In the case when both $X$ and $Y$ are affine schemes,
we have also the functor $\cxao\to C_Y^{old}(A')$. Having said this,
we possess the following collection of functors between categories of
sheaves on $X$ and on $Y$.
\roster
\item $f_{\bullet}:Sh_Y(A')\to\shx$.
\item The left adjoint of $f_{\bullet}$: $\fd:\shx\to Sh_Y(A')$.
This functor is exact if $f$ is flat.
\item (for $X$ and $Y$ affine) $\fs:Sh_Y^{old}(A')\to\shxo$.
\item (also for $X$ and $Y$ affine) The left adjoint of the
previous functor, denoted by $f^{*}$. This functor is also
exact if $f$ is flat.
\endroster

\sbh{2.3} Our next goal will be to define a certain subcategory $\qcsh$
in $\shx$ which we will call the category of $\qc$ sheaves.
The category $\qcsh$ will have properties analogous to those of the category
of $\qc$ sheaves of $O(X)$-modules inside the category of all sheaves of
$O(X)$-modules over a scheme $X$.
It will turn out that for an affine scheme, the category
of $\qc$ sheaves is equivalent to $\shxo$.

\sbh{2.3.1} Let $X$ be an affine scheme.
We have the natural inclusion functor $\cxao\to\cxan:B\to(X,B)$.
In this case there is the direct image
functor denoted $(\no):\shxn\to\shxo:\Gamma(B,(\no)(S))=\Gamma((X,B),S)$,
and the left adjoint of  $(\no)$, denoted by $(\on)$. The
functor $(\on)$ is exact and the functor $(\no)$ is left exact.
$R^{\bullet}(\no)$ will denote the right derived functor of $(\no)$.

\sbh{2.3.2}
Let now $f$ be a map $(Y,A')\to(X,A)$ with $X$ and $Y$ affine.
\proclaim{Lemma}
\roster
\item
The functors
$\fs\circ(\no)\text{ and } (\no)\circ f_{\bullet}:Sh_Y^{new}(A')\to\shxo$
are canonically
isomorphic.
\item
The same for the functors
$(\on)\circ f^{*}$ and $\fd\circ(\on)$ from $\shxo$ to
$Sh_Y^{new}(A')$.
\endroster
\endproclaim

The proof is obvious.

\sbh{2.3.3}
We will now describe the functor $(\on)$ more explicitly.
The next result can be considered as an analog of Serre's lemma.
\proclaim{Theorem}

Let $X$ be an affine scheme. Then the adjunction morphism of functors
$Id_{\shxo}\to(\no)\circ(\on)$ is an isomorphism.
\endproclaim
\demo{Proof of the Theorem}

Let $X$ be an affine scheme and consider a full subcategory ${\cxan}_{aff}$ of
\newline
$\cxan$ formed by pairs $(U,B_U)$ with $U$ affine. This subcategory
carries a natural Grothendieck topology. Let $emb$ be the embedding
functor $emb:{\cxan}_{aff}\to\cxan$. It is clearly a functor between sites
in the sense of 1.1.2.

The following lemma is easy to prove.

\proclaim{Lemma}
The functor $emb_{\bullet}:Sh(\cxan)\to Sh({\cxan}_{aff})$ is an equivalence
of categories. In particular,
$$emb_{\bullet}\circ emb^{\bullet}\simeq Id \text{  and  }
emb^{\bullet}\circ emb_{\bullet}\simeq Id $$
\endproclaim

Let us start now with a sheaf $S\in\shxo$ and consider the following
presheaf $S'$ on ${\cxan}_{aff}$: for $j:(U,B_U)\to (X,A)$ we set
$$\Gamma((U,B_U),S')\simeq\Gamma((U,B_U),j^{*}(S)).$$

We claim that this presheaf is in fact a sheaf canonically isomorphic
to $emb_{\bullet}(\on)(S)$. Indeed, without restricting
generality it suffices to check that if $j:(U,B_U)\to (X,A)$ is a covering,
the complex
$$0\rightarrow\Gamma((X,A),S'))\to\Gamma((U,B_U),S')
\to\Gamma((U_1,{B_U}_1),S')$$
is exact at first two places, where $j_1:(U_1,{B_U}_1)\to (X,A)$ denotes
the fiber product of $(U,B_U)$ with itself over $(X,A)$.

However, we know by 1.4.1, that
the complex of sheaves in $\shxo$
$$0\rightarrow S\to j_{*}j^{*}(S)\to {j_1}_{*}j_1*(S)$$
is exact at first two places.
The complex of groups, whose exactness we are proving, is just
the complex of global sections of this complex of sheaves and it is
exact at first two places  because the functor of taking sections is left
exact.

In order to establish the isomorphism $S'\to emb_{\bullet}(\on)(S)$,
we must exhibit an isomorphism:
$$Hom_{Sh({\cxan}_{aff})}(S',emb_{\bullet}(M))\to
Hom_{\shxo}(S,(\on)(M))$$ for any $M\in\shxn$. The latter is,
however, clear from the construction of $S'$.

In order to finish the  proof of the theorem it remains to show that
$(\no)\circ emb^{\bullet}(S)'\simeq S$ but this is obvious.

\enddemo

\sbh{2.3.4}

Let us now present several corollaries of the above theorem.

\proclaim{Corollary 1}
Let $X$ be an affine scheme and let $S\in\shxo$.
Let also $j:(U,B_U)\to (X,A)$ be an object of $\cxan$
with $U$ affine. Then
$$\Gamma((U,B_U),(\on)(S))\simeq\Gamma((U,B_U),j^{*}(S))$$
\endproclaim

\demo{Proof of Corollary 1}
We have
$$\align
\Gamma((U,B_U),(\on)(S))\simeq\Gamma((U,B_U),(\on)(S|(U,B_U)))\simeq  \\
\Gamma((U,B_U),j^{\bullet}(\on)(S))\simeq\Gamma((U,B_U),(\on)(j^{*}(S)))\simeq
\\
\Gamma((U,B_U),(\no)(\on)(j^{*}(S)))\simeq\Gamma(B_U,j^{*}(S))
\endalign$$
Here the first two isomorphisms follow from 1.3.1 Example 4, the third one
is a consequence of Lemma 2.3.2, the fourth isomorphism is the definition
of the functor $(\no)$ and the last one follows from the Theorem.
\enddemo{QED}

\proclaim{Corollary 2}
The functor $(\on)$ realizes $\shxo$ as a full abelian subcategory
of $\shxn$ stable under extensions.
\endproclaim
This is a formal consequence of the Theorem.

\proclaim{Corollary 3}
For any map $f:(Y,A')\to(X,A)$ with $Y$ and $X$ affine and for
any $S'\in\shxn$ the canonical morphism
$$f^{*}\circ(\no)(S')\to(\no)\circ\fd(S')$$
is an isomorphism provided that $S'$ is isomorphic to $(\on)(S)$
for some $S\in\shxn$.
\endproclaim
\demo{Proof of Corollary 3}

The canonical natural transformation:
$$f^{*}\circ(\no)\to(\no)\circ\fd:$$  follows from the standard mutual
adjunction properties of the functors $\fd$, $f_{\bullet}$, $f^{*}$, $\fs$,
$\no$,
$\on$.

If we put now $S'=(\on)(S)$, we will  get on the left hand side
$f^{*}\circ(\no)(\on)(S)\simeq f^{*}(S)$ whereas on the right hand side we get
$$(\no)\circ\fd\circ(\on)(S)\simeq (\on)\circ(\on)\circ f^{*}(S)\simeq
f^{*}(S)$$
and it is easy to verify, that under these identifications the above
natural transformation yields the identity morphism on $f^{*}(S)$.
\enddemo{QED}

\proclaim{Corollary 4}

The converse of the previous corollary is true:
if $X$ is an affine scheme
and if $S\in\shxn$, then $S\in\shxo$ if and only if
for every pair $(Y,A')$ mapping to $(X,A')$ by means of
$f$ with $Y$ being also affine,
the canonical map $f^{*}\circ(\no)(S)\to(\no)\circ\fd(S)$ is an isomorphism.
\endproclaim
\demo{Proof of Corollary 4}
Let $(\on)(\no)(S)\to S$ be the canonical adjunction morphism.
It follows then from Corollary 1 and from the assumption, that
for any $(U,B_U)\in\cxan$ with $U$ affine the above map induces
an isomorphism $\Gamma((U,B_U),(\on)(\no)(S))\simeq\Gamma((U,B_U),S)$.

Hence $(\on)(\no)(S)\simeq S$ and the assertion follows.
\enddemo{QED}

\proclaim{Corollary 5} Let $X=Spec(R)$ be an affine scheme.
\roster
\item
$R^i(\no)(\on)(S)=0$ for any $S\in\shxo$ and for any $i\geq 1$.
\item
The functor
$(\on)=\P:D(\shxo)\to D(\shxn)$ if fully faithful.
\endroster
\endproclaim
\demo{Proof of Corollary 5}

To prove the first point,we must show that if
$$0\to S\to K_1\to K_2\to 0$$
is an exact sequence of sheaves in $\shxn$ with $S\in\shxo$,
the sequence
$$0\to S\to (\no(K_1))\to (\no(K_2))\to 0$$
is exact in $\shxo$.

For this it suffices to check that if $B\in\cxao$ is a free algebra
built on a projective $R$-module, the complex of sections
$$0\to\Gamma(B,S)\to\Gamma(B,K_1)\to\Gamma(B,K_2)\to 0$$
is exact.

Now, since $\shxo\in\shxn$ is stable under extensions,

$Ext^1_{\shxo}(Const_{(X,B)},S)=0$ implies $H^1_{\shxn}((X,B),S)=0$
and the assertion follows.

The second point readily follows from the first one.

\enddemo{QED}
\sbh{2.3.5}
It is worthwhile to notice the following fact:
\proclaim{Lemma} Let $X$ be as in the Proposition.
\roster
\item
A sheaf $S\in\shxn$ comes from $\shxo$ if and only if the canonical
map $(\on)\circ(\no)(S)\to S$ is an isomorphism.
\item A sheaf $S\in\shxn$ comes from $\shxo$ if and only if for
any/some pair $(Y,A')\cxao$ covering $(X,A)$, $\fd(S)$
comes from $Sh_Y^{old}(A')$.
\endroster
\endproclaim
\demo{Proof}
The first point is clear and it implies the second one
because of 2.3.2 and corollary 3.
\enddemo{QED}

\sbh{2.4}
We arrive now to the definition of a $\qc$ sheaf in $\shx$.
\sbh{Definition}
Let $X$ be first an affine scheme. A sheaf $S\in\shxn$ is said to be
$\qc$ if it belongs to $\shxo$.

Let $X$ now be an arbitrary scheme. A sheaf $S\in\shx$ is said to be
quasi-coherent if for some (in fact, for any, by Lemma 2.3.5)
pair $(U,B_U)\in\cxan$ covering $(X,A)$ with $U$ affine, the restriction
of $S$ onto $(U,B_U)$ is $\qc$ in the sense of the previous definition.

Quasi-coherent sheaves form a full abelian subcategory in $\shx$,
which will be denoted by $\qcsh$. For $X$ affine, this category
coincides with $\shxo$.
\sbh{2.4.1 Example}
Let $(X,B)\in\cxa$, then $Const_{(X,B)}\in\qcsh$.

\sbh{2.4.2}
The category $\qcsh$ has enough injective objects because every
$Sh^{qc}_U(B_U)$ with $U$ affine does (cf. [Ar,Gr]).
\sbh{2.5}
We have the natural functor $\P:\Ds\to D(\shx)$  that sends \newline
$\Ds$
to the subcategory $\Dsq$ that consists of objects of $D(\shx)$
with quasi-coherent cohomologies.
\proclaim{Theorem}
The above functor induces an equivalence of categories: \newline
$\P:\Dbs\to\Dbsq$ \endproclaim

Corollary 5 of Theorem 2.3.3 implies the assertion for $X$ affine,
as well as the following lemma:
\proclaim{Lemma}
Let $f:(Y,A')\to (X,A)$ be a map such that
\roster
\item
$f:Y\to X$ is an affine morphism of schemes.
\item
$A'\to f^*(A)$ is an isomorphism.
\endroster

Then $R^if_{\bullet}(S)=0$ for any $i\geq 1$ and for any
$S\in Sh^{qc}_Y(A')$.
\endproclaim
\demo{Proof of the Lemma}

$R^if_{\bullet}(S)$ is a sheaf associated to the presheaf
$$(U,B_U)\to H^i((f^{-1}(U),f^*(A)|(f^{-1}(U)),S).$$
Now, for every $U$ each $B_U$ can be covered by a free one and
Corollary 5 of 2.3.3 finishes the proof.
\enddemo{QED}

\demo{Proof of the Theorem}

Since any complex is glued from its cohomologies, it suffices to prove that
for $S_1,S_2\in\qcsh$ the map $\P:Ext^i(S_1,S_2)\to Ext^i(\P(S_1),\P(S_2))$
is an isomorphism.

Choose $j:U\to X$ to be a covering in $Zar_X$ with $U$ affine.
Then $j:U\to X$ is an affine morphism since $X$ is a separated
scheme. Choose also an embedding $j^{\bullet}(S_2)\to I$ where $I$ is an
injective
object of $Sh^{qc}_Y(j^{*}(A))$.

$S\to j_{\bullet}(I)$ is an injection and let $K$ denote the cokernel. By
the above Lemma,
$$R^{\bullet}j_{\bullet}(\P(I))\simeq\P(R^{\bullet}j_{\bullet}(I))\simeq
\P(j_{\bullet}(I))$$
we have a commutative diagram
$$
\CD
Ext^{i-1}(S_1,K)  @>{\sim}>>  Ext^i(S_1,S_2)        \\
@V{\P}VV                             @V{\P}VV        \\
Ext^{i-1}(\P(S_1),\P(K))  @>{\sim}>>  Ext^i(\P(S_1),\P(S_2))
\endCD
$$
and the assertion follows by induction on $i$.
\enddemo{QED}

\heading 3. $A$-bimodules and sheaves on $\cxa$ \endheading
\sbh{3.0}
In this section we will study the connection of the category of
$\qc$ sheaves of $A$-bimodules to that of $\qc$ sheaves on $\cxa$.
The material here is parallel to the one of Section 3 in [Ga]. The category
of $\qc$ sheaves of $A$-bimodules will be denoted by $\Aq$.
\sbh{3.1}
Let us recall several definitions from [Ga].
If $B$ is a $\qc$ sheaf of algebras on a scheme $X$, we denote by $I_B$ the
sheaf of $B$-bimodules given by $I_B=ker(B\otimes B\to B)$ (the map here is
the multiplication).

If $M$ is a $\qc$ sheaf of $B$-bimodules, we will denote by $\Omega(B,M)$
the group $Hom(I_B,M)$.

\sbh{3.2} Let now $X$ be a scheme and let $A$ be a $\qc$ sheaf of
associative algebras on $X$. We will construct a localization functor
$\Im:\Aq\to\qcsh$:

Let $M\in\Aq$. Consider the presheaf $\Im(M)$ on $\cxa$ given by
$$\Gamma((U,B_U),\fm)=\Omega(B_U,M|U).$$

The following Lemma is proven by a straightforward verification.

\proclaim{Lemma} This presheaf is in fact
a sheaf.
\endproclaim

\sbh{3.2.1}
If $X$ is an affine scheme, similar constructions can be carried out
in the category $\cxao$. In this case we denote the localization functor
by ${\Im}^{old}$.
\proclaim{Lemma}
The functors $$(\no)\circ\Im\text{ and }{\Im}^{old}:\Aq\to\qcsh$$
are canonically isomorphic.
\endproclaim
The proof is obvious.

\sbh{3.3}
Let us describe $\fm$ in a slightly different way. Consider the sheaf of
algebras $A\oplus M$ over $X$, $(X,A\oplus M)\in\cxa$. Then it is a group-like
object
in this category: $$Hom((U,B_U),(X,A\oplus M))=\Omega(B_U,M|U),$$
and $\Im(M)$ is a sheaf given by $\Gamma((U,B_U),\Im(M)))
=Hom((U,B_U),(X,A\oplus M))$.
In other words, $\Im(M)$ is a group like object in the category of
sheaves of sets with $\Im(M)=Const'_{A\oplus M}$. (Warning:
note the difference between $Const'_X$ and $Const_X$: the former is the
constant
sheaf corresponding to $X$ in the category of sheaves of sets, whereas the
latter
is the constant sheaf in the category of the sheaves of groups, cf. 1.2 and
1.3.1(2).)

\sbh{3.3.1}

\proclaim{Proposition}
\roster
\item
Let $f:(Y,A')\to (X,A)$ be a map with $Y$ and $X$ affine. Then the functors
$${\Im}^{old}\circ f^{*}\text{ and }f^{*}\circ {\Im}^{old}:\Aq(X)\to
Sh^{old}_{A'}(Y)$$
are canonically isomorphic.
\item
For $X$ and $Y$ arbitrary and for any ,$f:(Y,A')\to (X,A)$ the functors
$${\Im}\circ f^{*}\text{ and }f^{\bullet}\circ {\Im}:\Aq(X)\to
Sh_{A'}(Y)$$
are canonically isomorphic.
\item
Let once again $X$ be affine. Then the functors
$$(\on)\circ {\Im}^{old}\text{ and }\Im:\Aq\to Sh^{new}_A(X):\Aq(X)\to
Sh_{A}(X)$$
are canonically isomorphic.
\item
For $X$ arbitrary the functor $\Im:\Aq\to Sh_A(X)$ takes values in $\qcsh$.
\item
The functor $\Im:\Aq\to\qcsh$ is exact and faithful.
\endroster
\endproclaim

\demo{Proof of the Proposition}
First two points are immediately deduced from the following general Lemma:

\proclaim{Lemma 1}
Let $F:C_1\to C_2$ be a functor between two sites. Let $A\in C_1$
be an (abelian) group-like object and let also $Const_A$ be the sheaf of
abelian groups associated with it. Suppose that $F(A)$ is an (abelian)
group-like object in $C_2$ as well and that $F:End(A)\to End(F(A))$
is a homomorphism of groups. Then $F^{\cdot}(Const'_A)\simeq Const'_{F(A)}$.
\endproclaim
\demo{Proof of the Lemma 1} The proof follows from the following observation:
for any $S\in Sh(C)$, $$Hom(Const'_A,S)=\{\gamma\in\Gamma(A,S)|n\cdot\gamma=
n^*(\gamma)\}\text{ for any }n\in{\Bbb Z},$$
where $n$ on the right hand side denotes the endomorphism $n\cdot Id_A\in
End(A)$.
\enddemo{QED}

The third point of the Proposition follows from the first two and from
the Corollary 1 of 2.3.3.

(4) is an immediate consequence of (2) and of (3). In order to prove (5)
we may assume $X$ to be affine, where the assertion follows from the following
lemma, whose proof is a straightforward verification.
\proclaim{Lemma 2}
Let $X=Spec(R)$ and let $V$ be an $R$-module. Let also $Free_R(V)\in\cxao$
be the free associative algebra built on $V$. Then we have a canonical
isomorphism of functors:$\Aq\to Ab$
$$M\to Hom_R(V,M) \text{ and } M\to\Gamma(Free_R(V),\Im(M)).$$
\endproclaim
\enddemo{QED}

\sbh{3.3.2}
The following assertion is easy:
\proclaim{Lemma}
Let $f$ be a map from a pair $(Y,A'\simeq f^*(A))$ to the pair $(X,A)$.
We have then the direct image functor $f_*:A'_{qc}-mod\to\Aq$.
The functors $\Im_X\circ f_*$ and $f_{\bullet}\circ\Im_Y$ are canonically
isomorphic. \endproclaim

\proclaim{3.4 Lemma-Definition}
The functor $\Im:\Aq\to\qcsh$ admits a left adjoint
denoted by $\LL:\qcsh\to\Aq$ \endproclaim
\demo{Proof}
For any $S\in\qcsh$ we must construct a $\qc$ $A$-module $\LL(S)$,
satisfying $Hom(\LL(S),M)\simeq Hom(S,\fm)$ functorially in $M$.
Because of 3.3.1, it will suffice to construct $\LL(S)$ locally
on $X$. This reduces us to the situation when $X$ is affine. In this
case, any object in $\qcsh$ is a quotient of a one of the type
$Const_{Free(V)}$. However, for $S=Const_{Free(V)}$, Lemma 2 of 3.3.1 implies
$\LL(S)=F(V)$, where
$F(V)$ denotes the free $A$-module built on $V$. We finish the proof
by quoting the following general (and trivial)
\proclaim{Sub-lemma}
Let $C_1$ and $C_2$ be two abelian categories and let $F:C_1\to C_2$
be an additive left exact functor between them. Suppose that $F$ admits
a partially defined left adjoint functor which is however defined
on a large collection of objects in $C_2$ (i.e. any object in $C_2$
is a quotient of a one from this collection). Then this left
adjoint is defined on the whole of $C_2$. \endproclaim
\enddemo{QED}
\sbh{3.4.1}
{}From 3.3.2 and the above Lemma-definition we deduce by adjunction the
following result:
\proclaim{Lemma}
Let $f$ be a map from a pair $(Y,A'\simeq f^*(A))$ to the pair $(X,A)$.
The functors $${\LL}_Y\circ f^{\bullet}\text{ and }f^*\circ{\LL}_X:
Sh^{qc}_A(X)\to A'_{qc}-mod(Y)$$ are canonically isomorphic.
\endproclaim
\sbh{3.4.2}
Since the functor $\Im$ is exact, it can be prolonged to
a functor between the corresponding derived categories:
$\Dq\overset def\to\simeq D(\Aq)\to\Ds$ which will be also denoted by $\Im$.
Our next aim is to show that the functor $\LL$ (which is obviously
right exact) can also be derived into the functor
$L^{\cdot}\LL:\Dms\to\Dmq$, which will be the left adjoint functor to
$\Im:\Dmq\to\Dms$.
When $X$ is affine, the argument of the above lemma proves everything,
since the sheaves $Const_{Free(V)}$ with $V$ being a projective
$O(X)$-module form a set of projective generators of $\qcsh$. However, in order
to treat the general case an additional argument is needed, since
objects of the derived category cannot be reconstructed just from the local
information.

\proclaim{3.4.3 Theorem}
Let $X$ be an arbitrary scheme and $A$ be a $\qc$ sheaf of algebras on $A$.
Then the functor $\LL$ can be derived into the functor
$L^{\bullet}\LL:\Dms\to\Dmq$
which moreover satisfies $$Hom(L^{\bullet}\LL(S^{\bullet}),M^{\bullet})\simeq
Hom(S^{\bullet},\Im(M^{\bullet}))$$ functorially in $S^{\bullet}\in\Dms$ and in
$M^{\bullet}\in D^+_{qc}(A)$. \endproclaim
\sbh{Remark}
The category $\qcsh$ is lacking objects that would be acyclic for the
functor $\LL$. The situation here is similar to that in [Bo,Be], when
one wants to define the direct image functor for $D$-modules. As in
[Bo], there are at least two ways to overcome this difficulty:
a more straightforward one is to go beyond the category $\qcsh$
and work with arbitrary sheaves. In this case there are enough acyclic
objects for the functor $\LL$, but the drawback of this approach is
that we will have to rely on the equivalence of the categories $\Dq$
and $D(A-mod)$ with $\qc$ cohomologies as well as on Theorem 2.6.
Another way is the one described below:

\demo{Proof of the Theorem}
Consider first a pair $(Y,A')$ with $Y$ is affine and let $S$ be a
quasi-coherent sheaf
on $C_Y(A')$. We will construct a canonical sheaf $Can'(S)$ mapping
surjectively
onto $S$ with $Can'(S)$ being acyclic for the functor $\LL$. Namely,
$Can(S)=\oplus\Gamma(B,S)\otimes Const_B$, the sum being taken over
isomorphism classes of objects in $\cxao$ with $B$ a free algebra on
a projective $O-$module. This construction has the following two properties:

{\bf 1)}
For any map of sheaves $S\to S'$ there is a canonical map $Can'(S)\to
Can'(S')$.

{\bf 2)}
If $f:(Z,f^*(A))\to (Y,A)$ is a morphism of pairs with $Y\text{ and }Z$ affine,
there exists a canonical map $f^{\bullet}(Can'(S))\to Can'(f^{\bullet}(S))$.

Thus any complex $S^{\bullet}$ of sheaves bounded from above
in $\qcsh$ admits a admits a {\it canonical} resolution $Can(S^{\bullet})$ by a
complex consisting of sheaves acyclic with respect to the functor $\Im$.

Let now $S^{\bullet}$ be a complex bounded from above on $X$ giving rise
to an object of $\Dms$, and choose $j:U\to X$ to be a covering in $Zar_X$ with
$U$
affine. Put $A'\simeq j^*(A)$ and
let also $U_i$ be as in 1.4.1. All these schemes are affine since $X$ is
assumed to be separated.

For each $U_i$, fix the canonical resolution $Can(S^{\bullet}|U_i)$
of $S^{\bullet}|U_i$ as above.

Then for each $i$ we can form a complex $\LL(Can(j_i^{\bullet}(S^{\bullet})))$
of quasi-coherent sheaves of $A|U_i$-bimodules on $U_i$.

For each $i,j$ we have a boundary map ${p^{\bullet}_j}:U_i\to U_{i-1}$

and therefore we have a map of complexes:
$${p^{\bullet}_j}^i\LL(Can(j_{i-1}^{\bullet}(S^{\bullet})))\to
\LL(Can(j_i^{\bullet}(S^{\bullet})))$$
which is easily seen to be a quasi-isomorphism by 3.4.1, since the functor
$${p^{\bullet}_j}^{i}:Sh^{qc}_{A'_{i-1}}(U_{i-1})\to Sh^{qc}_{A'_i}(U_i)$$
is exact.

We have a complex of complexes:
$$0\rightarrow j_*\LL(Can(j^{\bullet}(S)))\rightarrow
{j_1}_*\LL(Can(j_1^{\bullet}(S)))
\rightarrow {j_2}_*(Can(j_2^{\bullet}(S)))\rightarrow ...$$
or, in other words, a double complex ${\LL'}^{\bullet\bullet}(S^{\bullet})$
in $\Aq$, whose associated complex we denote by
$Ass({\LL'}^{\bullet\bullet}(S^{\bullet}))$.

It is now a standard exercise to check, that for each $i$ the canonical map of
complexes $$\LL(Can(j^{\bullet}(S^{\bullet})))\to
Ass({\LL'}^{\bullet\bullet}(S^{\bullet}))|U$$
is a quasi-isomorphism. This in turn implies that the functor
$$S^{\bullet}\to Ass({\LL'}^{\bullet\bullet}(S^{\bullet}))$$ is a well-defined
functor
$\Dms\to\Dmq$.

To prove that this functor is in fact $L^{\bullet}\LL$ we are looking for,
it remains to verify the adjunction property:
$$Hom(L^{\bullet}\LL(S^{\bullet}),M^{\bullet})\simeq
Hom(S^{\bullet},\Im(M^{\bullet}))\eqno{(*)}$$
For this we must construct the adjunction morphisms
$$S^{\bullet}\to\Im\circ Ass({\LL'}^{\bullet\bullet}(S^{\bullet}))$$
                             and
$$Ass({\LL'}^{\bullet\bullet}(\Im(M^{\bullet})))\to M^{\bullet}$$
This is done in the following way:
$$
\CD
Ass[{j_i}_{\bullet}j_i^{\bullet}(Can(S^{\bullet}))] @<{\sim}<<
Can(S^{\bullet})  @>{\sim}>> S^{\bullet} \\
@VVV  \\
Ass[{j_i}_{\bullet}Can(j_i^{\bullet} S^{\bullet})] \\
@VVV  \\
Ass[{j_i}_{\bullet}\Im\LL(Can(j_i^{\bullet}S^{\bullet}))] @>{\sim}>>
\Im\circ Ass[{j_i}_*\LL(Can(j_i^{\bullet} S^{\bullet}))]
\endCD$$
for the first adjunction map, and
$$
\CD
Ass[{j_i}_*\LL(j_i^{\bullet}(\Im(M^{\bullet})))] @<<<
Ass[{j_i}_*\LL\circ Can(j_i^{\bullet}(\Im(M^{\bullet})))] \\
@VVV  \\
Ass[{j_i}_*j_i^*(\LL\circ\Im(M^{\bullet}))] @>>>
Ass[{j_i}_*j_i^*(M^{\bullet})] @>{\sim}>> M^{\bullet}
\endCD
$$
for the second one.
It is now easy to verify, that the adjunction maps constructed above
give rise to $(*)$.
\enddemo{QED}
\sbh{3.4.4}
The following statement readily follows from Lemmas 3.3.2 and 2.2.5
by adjunction and is implicit in the Theorem:
\proclaim{Lemma}
Let $(Y,A')\to (X,A)$ be a map such that $f:Y\to X$ is affine and flat
and such that $A'\to f^*(A)$ is an isomorphism.
Then $${L^{\bullet}\LL}_Y\circ f^{\bullet}\simeq f^*\circ
{L^{\bullet}\LL}_X$$ as functors ${D^-_{qc}}_X(A)\to {D^-_{qc}}_Y(A')$.
\endproclaim
\sbh{3.5 Definition}
Let $(X,A)$ be as before: a scheme with a $\qc$ sheaf of algebras on it.
We define $T^{\bullet}(A)$ to be the object of $\Dq$ given by $L^{\bullet}\LL
(Const_A)$.

{}From the fact that $\LL$ is left exact we infer that $H^i(T^{\bullet}(A))$
vanishes for $i>0$ and that $H^0(T^{\bullet}(A))=I_A$.

$T^{\bullet}(A)$ will be called the cotangent complex of $A$.
If $M$ is a $\qc$ sheaf of $A$-bimodules, we denote by $H^i_A(M)$ the
groups $Ext^i(T^{\bullet}(A),M)$.

\sbh{3.5.1 Example}
Suppose that $A$ is flat over $O(X)$. It follows from the results of
Quillen [Qu], that $T^{\bullet}(A)\simeq I_A$.
Indeed, this is true for $X$ affine, and then we apply 3.4.4.

\heading{4. deformation theory}\endheading
\sbh{4.0}
This section is almost a word by word repetition of [Ga],
after we adopt certain modifications connected with the fact that
we are working over a scheme.
\sbh{4.1.0}
For a scheme $X$, $O_i(X)$ will denote the sheaf $O[t]/t^{i+1}\cdot O(X)$.
\sbh{4.1.1}
Let $A$ be a $\qc$ sheaf of associative algebras on $X$.
\sbh{Definition}

The category $Deform^i(A)$ is defined to have
as objects quasi-coherent sheaves of associative $O_i(X)$-algebras $A_i$,
endowed with an isomorphism $A_i/A_i\cdot t\simeq A$ such that
$Tor^{O_i(X)}_1(A_i,O(X))=0$. (In other words, we need that
$$ker(t:A_i\to A_i)=im(t^i:A_i\to A_i)\text{ identifies under
a natural map with } A).$$  Morphisms in this category are just $O_i(X)-$
algebras homomorphisms respecting the identifications with $A$ modulo $t$.
This category is obviously a groupoid. It is called the category of $i$-th
level deformations of $A$.

We have natural functors $Deform^{i+1}(A)\to Deform^i(A)$ given by
reduction modulo $t^{i+1}$. If $A_i$ is an object in $Deform^i(A)$,
we denote by $Deform^{i+1}_{A_i}(A)$ the category-fiber of the above functor.
This category, which is obviously a groupoid too, is called the category
of prolongations of $A_i$ onto the $i+1$-st level.
\sbh{4.1.2}
We are now ready to state the main result of the present paper:
\proclaim{Theorem}
\roster
\item
The category $Deform^1(A)$ is equivalent to the category $T_{\Im(A)}$
of $\Im(A)$-torsors
on $\cxa$.
\item
If $A_i\in Deform^i(A)$, and if $A_{i+1}\in Deform^{i+1}_{A_i}(A)$,
$Aut(A_{i+1})$ is canonically isomorphic to
$\Omega(A,A)=\Gamma(Const_{(X,A)},\Im(A))$.
\item
To any $A_i\in Deform^i(A)$ one can associate a gerbe $G_{A_i}$
bound by $\Im(A)$ on $\cxa$ in such a way that $G_{A_i}((X,A))$ is canonically
equivalent to $Deform^{i+1}_{A_i}(A)$.
\endroster
\endproclaim
\sbh{Remarks}
It is easy to notice that the first two points in the statement
of the Theorem are special cases of the third one.

\demo{Proof}
\roster
\item
The functor $T:Deform^1(A)\to T(\Im(A))$ is given by:
$$\Gamma((U,B_U),T(A_1))=O(U)-\text{algebras homomorhisms }B_U\to A_1|U$$
and it is easy to show that it is an equivalence
of categories.
\item
This is a direct verification.
\item
We define the gerbe $G_{A_i}$ as follows:

$G_{A_i}((U,B_U))$ is the groupoid of $O_i(U)$-algebras ${B_U}_{i+1}$
with an isomorphism  ${B_U}_{i+1}/t^{i+1}\cdot {B_U}_{i+1}\simeq
B_U\underset{A|U}\to\times (A_i)|U$ such that \newline
$ker(t^{i+1}:{B_U}_{i+1}\to {B_U}_{i+1})=im(t:{B_U}_{i+1}\to {B_U}_{i+1})$ and
identifies under a natural morphism with $(A_i)|U$.

Functors $G_{A_i}((U,B_U))\to G_{A_i}((V,C_V))$ for maps $(V,C_V)\to (U,B_U)$
are given by taking fiber products.

The fact that $G_{A_i}$ is indeed a gerbe bound by
$\Im(A)$ over $\cxa$ and that its fiber over $(X,A)$ is equivalent to
$Deform^{i+1}_{A_i}(A)$ can be easily verified.
\endroster
\enddemo{QED}
\sbh{4.2}
We will now translate the assertions of the above theorem into cohomological
terms.
\sbh{Remark}
It has been proven (Theorem 2.5) that cohomologies of a quasi-coherent sheaf
over $\cxa$ computed either in the category of all sheaves or in the
category of quasi-coherent sheaves coincide. Hence classes of gerbes and
torsors
over  $\qc$ sheaves can be computed inside the $\qc$ category.

\sbh{1-st Level Deformations}
The groupoid $Deform^1(A)$,
which is a priori endowed with a monoidal structure is a Picard category
with $\pi_0(Deform^1(A))$ (the set of isom. classes of objects)
being a group isomorphic to $Ext^1(Const_{(X,A)},\Im(A))\simeq H^1_A(A)$.
This follows from 1.5.2, from 3.4 and from 4.1.2(1).
\sbh{Prolongation of Deformations 1}
If $A_i$ is an $i$-th level deformation, there exists a canonical class
in $H^2_A(A)$ which is zero if and only if there exists a prolongation
$A_{i+1}$ of $A_i$. This follows from 1.5.5, from 3.4, and from 4.1.2(3).
\sbh{Prolongation of Deformations 2}
Suppose that for a given $i$-th level deformation $A_i$ the category
$Deform^{i+1}_{A_i}(A)$ has an object. Then $\pi_0$ of this category
is a torsor over the abelian group $H^1_A(A)$. This follows from 1.5.1,
from 3.4. and from 4.1.2(3).

\sbh{4.3 Example}
Suppose now that the sheaf $A$ is flat as a sheaf of $O(X)-$modules.
{}From Example 3.5.1, it follows that deformations of $A$ are controlled
by $Ext^i(I_A,A)$ ($Ext$s being taken in the category of $\qc$ sheaves of
$A$-bimodules ) for $i=1,2$.

\enddocument